\documentclass[12pt]{amsart}

\usepackage{latexsym}
\usepackage{amsmath}
\usepackage{amssymb}
\usepackage{amscd}
\usepackage{psfig}
\usepackage{multicol}

\newtheorem{theorem}{Theorem}[section]

\newtheorem{lemma}{Lemma}[section]

\def\Q{\mathbb{Q}}

\def\R{\mathbb{R}}
\def\C{\mathbb{C}}
\def\iso{\cong}

\def\t{\mathfrak{t}}
\def\g{\mathfrak{g}}
 \def\k{\mathfrak{k}}
 \def\cO{\mathcal{O}}

\title[Symplectic reductions]{Cohomology of symplectic reductions of generic coadjoint orbits}

\date{ \today }

\author[Goldin]{R.\ F.\ Goldin}
\address[R.\ F.\ Goldin]{Mathematical Sciences\\
George Mason University\\ MS 3F2, 4400 University Dr.\\ Fairfax, VA 22030}
\email{rgoldin@gmu.edu}

\author[Mare]{A.-L.\ Mare}
\address[A.-L.\ Mare]{Department of Mathematics\\ University of Toronto
 \\Toronto, Ontario M5S 3G3, Canada}
 \email{amare@math.toronto.edu}

\begin{document}
\begin{abstract}Let $\cO_\lambda$ be a generic coadjoint orbit of a compact semi-simple Lie group $K$. Weight varieties are the symplectic reductions of $\cO_\lambda$ by the maximal torus $T$ in $K$. We use a theorem of Tolman and Weitsman to compute the
cohomology ring of these varieties. Our formula relies on a {\em Schubert basis} of the equivariant cohomology of $\cO_\lambda$ and it makes explicit the dependence on $\lambda$ and a parameter in $Lie(T)^*=:\t^*$.
\end{abstract}
\maketitle

\section{Introduction}

Let $K$ be a compact semisimple Lie group, $T\subset K$ a maximal
torus and $\t \subset \k$ their Lie algebras. Pick a fundamental
chamber in $\t^*$ and a point $\lambda$ in the interior. Let
$\cO_\lambda$ be the  orbit of $\lambda$ under the coadjoint
representation
 of $K$ on $\k^*$. $\cO_\lambda$ is diffeomorphic to the flag variety $K/T$ and it has a naturally occurring symplectic form $\omega$ known as the Kirillov-Kostant-Souriau form.
The action
of $T$ on $\cO_{\lambda}$ is Hamiltonian, which means that there is an invariant map
$$\Phi : \cO_{\lambda}\to \t^*$$
satisfying $\omega(X_\eta,\cdot)=d \Phi^\eta$, where $\eta\in\t$, $X_\eta$ the vector field on $\cO_\lambda$ generated by $\eta$, and $\Phi^\eta(m)=\Phi(m)(\eta)$ defined by the natural pairing between $\t$ and $\t^*$. We call $\Phi$ a {\em moment map} for this action.

The image of $\Phi$ is the convex hull of $W\cdot\lambda$, the Weyl group
orbit of $\lambda$.
Let $\mu\in \Phi(\cO_{\lambda})$ be a regular value of $\Phi$. We define the {\em symplectic reduction} at $\mu$ by
$$\Phi^{-1}(\mu)/T=\cO_{\lambda}//T(\mu).$$

The goal of this note is to give a presentation of the
cohomology\footnote{Only cohomology with coefficients in the field
$\Q$ of rational numbers will be considered throughout this
paper.} ring of $\cO_{\lambda}//T(\mu)$ in terms of the root
system of $K$. We present $H^*(\cO_\lambda/\!/T(\mu))$ as a quotient of the $T$-equivariant cohomology ring $H_T^*(\cO_\lambda)$ by a certain ideal. We rely on the following fundamental result.

\begin{theorem}[Kirwan] Let $M$ be a compact symplectic manifold with a Hamiltonian $T$ action, where $T$ is a compact torus. If
$\mu\in \t^*$ is a regular value of $\Phi$,
 then the restriction map in equivariant cohomology
$$\kappa :H_T^*(M)\to H^*_T(\Phi^{-1}(\mu))$$
is surjective.
\end{theorem}As the $T$ action is locally free on level sets of the moment map at
regular values, $H_T^*(\Phi^{-1}(\mu))=H^*(M//T(\mu))$. The
resulting map $\kappa:H_T^*(M)\to H^*(M//T(\mu))$ is called the
{\em Kirwan map}.  Kirwan's result is of particular importance
because the equivariant cohomology can be described in terms of
the equivariant cohomology of the fixed point sets of the $T$
action. In the case of isolated fixed points, this is just a sum
of polynomial rings.

\begin{theorem}[Kirwan]\label{th:restriction} Let $M$ be a compact Hamiltonian $T$-space. Let $M^T$ denote the fixed point set of the $T$ action. The restriction map
$$i^*:H_T^*(M)\to H_T^*(M^T)$$
is injective. In the case that $M^T$ is a finite set of points,
$H_T^*(M^T)=\oplus_{p\in M^T} \Q[x_1,\dots,x_n]$ where $n=\dim T$.
\end{theorem}

A  presentation of the cohomology ring of the reduced space
$M//T(\mu)$ can be obtained by using the following  description of the kernel
of the Kirwan map,
which is due to Tolman and Weitsman \cite{To-We}.
If $\alpha$ is in $H_T^*(M)$ we denote
$${\text supp} (\alpha)=\{p\in M^T: \alpha|_p\neq 0\}$$
 Fix an arbitrary inner
product $\langle \ , \ \rangle$ on $\t^*$.

\begin{theorem}[Tolman-Weitsman] The kernel of the Kirwan
map $\kappa$ is the ideal of $H_T^*(M)$ generated by all
$\alpha \in H_T^*(M)$ with the property that there exists $\xi\in \t^*$
such that
$$\Phi({\text supp} (\alpha))\subset \{x\in \t^* | \langle \xi,x\rangle
\leq \langle \xi, \mu\rangle\}.$$
In other words, $\alpha$ is in $\ker \kappa$ if and only if all points of
${\text supp} (\alpha)$ are mapped by
 $\Phi$ to the same side of an affine hyperplane  in $\t^*$ which passes
through $\mu$.
\end{theorem}
The $T$-equivariant cohomology ring of the coadjoint orbit
$\cO_{\lambda}=K/T$ is well understood.  Kostant and Kumar
constructed in  a basis $\{x_w\}_{w\in W}$ of $H_T^*(K/T)$ as a
$H_T^*(pt)$-module, which we refer to as the {\it Schubert basis}
\cite{Ko-Ku}. Let $B$ be a Borel in $G:=K^\C$, and $B_-$ an
opposite Borel. For any $v\in W$, let
$X_v=\overline{B_-\tilde{v}B}/B$, where $\tilde{v}$ is any choice
of lift of $v\in W$ in the normalizer of the torus. These {\em
opposite Schubert varieties} are $T$-invariant subvarieties of
$G/B\iso K/T$. The basis $\{x_w\}$ is uniquely defined by the
property that
$$
\int_{X_v}x_w=\delta_{vw}.
$$
Theorem \ref{th:restriction} suggests the importance of knowing
how to restrict the classes $x_w$ to fixed points $W\cdot
\lambda$. This formula was worked out for general $K$  by S.
Billey \cite {Bi}. In particular, it is easy to show that
$x_w|_v=0$ if $v\not\leq w$ in the Bruhat order.\footnote{The
class $x_w$ differs from the $\xi^w$ constructed in \cite{Ko-Ku}
by the relationship $x_w:=w_0\cdot\xi^{w_0w}$, where $w_0$ is the
longest element of $W$.} In other words,
$$supp(x_w)=\{v\lambda : v\leq w\}.$$
To each $\tau \in W$ we can associate the new basis
$$\{x_w^{\tau}=\tau \cdot x_{\tau^{-1}w}\}_{w\in W},$$
whose elements have the property
 $$supp(x_w^{\tau})=\{v\lambda : \tau^{-1}v\leq \tau^{-1}w\}.$$

Let  $\lambda_1, \ldots , \lambda_l\in \t^*$ denote the fundamental
weights associated to the chosen fundamental chamber of $\t^*$. Let
 $\langle \ , \ \rangle$ be the restriction to
$\t^*$ of a $K$-invariant product on $\k^*$. Our main result is:

\begin{theorem}\label{th:main} The cohomology ring $H^*(\cO_{\lambda}//T(\mu))$ is
isomorphic to the
quotient of $H_T^*(K/T)$ by the ideal generated by
$$\{x_v^{\tau}: {\text there \  exists \ }  j  \ {\text such \ that } \
\langle \lambda_j,\tau^{-1}v\lambda\rangle \leq
\langle \lambda_j, \tau^{-1}\mu\rangle\}.$$
\end{theorem}

\noindent {\bf Remarks.}

\noindent 1. One can take the description of $H_T^*(K/T)$ (see for instance \cite{Br})
and deduce a  precise presentation of the cohomology ring $H^*(\cO_{\lambda}//T(\mu))$
in terms of generators and relations.

\noindent 2. For $K=SU(n)$ this result was proven by the first
author in \cite{Go1}.

\vspace{0.2cm}

\noindent {\bf Acknowledgement.} The second author would like to
thank Lisa Jeffrey for introducing him to the topic of the paper.
Both authors would like to thank her for a careful reading of the
manuscript and for suggesting several improvements.

\section{Primary description of $\ker \kappa$}
For any $\xi \in \t^*$ we denote by $f_{\xi}$ the corresponding
height function on $\cO_{\lambda}$,
$$f_{\xi}(x)=\langle \xi,x\rangle.$$
 Under the pairing between $\t^*$ and $\t$, the function $f_\xi$
is a component of the moment map. In fact, it is well known that
$f_\xi$ is Morse-Bott
 for all $\xi\in\t^*$.
Denote by $C \subset \t^*$  the fundamental (positive) Weyl
chamber, which can be described by
$$ C=\{r_1\lambda_1 +\ldots +r_l\lambda_l: {\text \ all } \ r_j >0\},$$
and let $\overline{C}$ be its closure.

\begin{lemma} Let $\tau$ be in $W$ and $\xi$ in $\tau\overline{C}$.
  If $\tau^{-1}v < \tau^{-1}w$ in the Bruhat order, then $f_{\xi}(v\lambda)\leq
f_{\xi}(w\lambda)$.
\end{lemma}

\begin{proof}  The result follows immediately from the fact that if
$\xi \in C$, then the unstable manifold of $f_{\xi}$ through
$v\lambda$ with respect to the K\"ahler metric on
$$\cO_{\lambda}=K/T=G/B$$ is just the Bruhat cell $B\cdot vB/B$ (see for instance \cite{Koch}).
\end{proof}

The main result of this section is:
\begin{theorem} Suppose that $\alpha\in H_T^*(\cO_{\lambda})$ has the
property that
$$\Phi({\text supp} (\alpha))\subset \{ x\in \t^*: \langle \xi, x\rangle
\leq \langle \xi, \mu\rangle \}.$$
Then $\alpha$ can be decomposed  as
$$\alpha  =\sum_{w\in W}a_w^{\tau}x_w^{\tau} $$
with $a_w^{\tau}\in H_T^*(pt)$, such that if $a_w^{\tau}\neq 0$
then
$$\Phi(supp (x_w^{\tau}))\subset  \{ x\in \t^*: \langle \xi, x\rangle
\leq \langle \xi, \mu\rangle \}.$$
\end{theorem}
\begin{proof}
 Take $\tau\in W$ such that $\xi\in
\tau\overline{C}$. Suppose that the decomposition of $\alpha$ with respect
to the basis $\{x_w^{\tau}\}_{w\in W}$ is of the form
\begin{equation}
\alpha=\sum_{w\in W}a_w^{\tau} x_w^{\tau}+a_{v_1}^{\tau}x_{v_1}^{\tau}+\ldots +
a_{v_r}^{\tau}x_{v_r}^{\tau},
\end{equation}
where the first sum contains only $w$ with
$$\langle \xi,w\lambda \rangle
\leq \langle \xi, \mu\rangle,$$
whereas
$$\langle \xi,v_j\lambda \rangle
> \langle \xi, \mu\rangle, \quad  a_{v_j}^{\tau}\in S(\t^*), a_{v_j}^{\tau}\neq 0,$$
for any $1\leq j \leq l$. We may assume that $v_1$ has the
property that there exists no $j> 1$ with $\tau^{-1}v_1
<\tau^{-1}v_j$. Now let us evaluate both sides of (1) at
$v_1\lambda$. Since
$$\langle \xi,w\lambda \rangle
\leq \langle \xi, \mu\rangle < \langle \xi,v_1\lambda \rangle,$$
by Lemma 2.1 we must have
$$x_w^{\tau}|_{v_1\lambda}=0$$
for any $w$ corresponding to a term in the first sum in (1). It follows that
$$\alpha|_{v_1\lambda}=a_{v_1}^{\tau}x_{v_1}^{\tau}|_{v_1\lambda}\neq 0$$
so $v_1\lambda$ is in $supp(\alpha)$ even though $\langle \xi,v_1\lambda \rangle
> \langle \xi, \mu\rangle$. This is a contradiction.
\end{proof}

\section{Proof of the main result}

We now prove Theorem~\ref{th:main}. Let $v$ and $\tau$  in
$W$ be such that
\begin{equation}
\langle \lambda_j ,\tau^{-1}v\lambda \rangle \leq
\langle \lambda_j ,\tau^{-1}\mu \rangle,
\end{equation}
for some $1\leq j \leq l$. We show that $x_v^{\tau}$ is in the kernel of the Kirwan
map $$\kappa:H_T^*(\cO_\lambda)\to H^*(\cO_{\lambda}/\!/T(\mu)).$$
Let $\xi=\tau \lambda_j$ be in $\tau\overline{C}$. Note that if
$w\in supp (x_v^{\tau})$, then $\tau^{-1}w \leq \tau^{-1}v$ implies by
 Lemma 2.1
$$\langle \xi,w\lambda \rangle \leq \langle \xi,v\lambda \rangle \leq \langle \xi, \mu\rangle.$$
Thus $x_v^\tau\in\ker\kappa$.

Now let us consider $\alpha\in H_T^*(K/T)$ with the property that there
exists $\xi\in \t^*$ with
$$supp(\alpha)\subset \{x\in \t^*|\langle \xi, x\rangle \leq
\langle \xi, \mu \rangle\}.$$
Take $\tau\in W$ such that $\xi \in\tau\overline{C}$.
By Theorem 2.2, we can decompose $\alpha$ as
\begin{equation}\label{eq:decompose}
\alpha=\sum_{w\in W}a_w^{\tau}x_w^{\tau}
\end{equation}
where  $a_w^{\tau}$ can be nonzero only if
$$supp(x_w^{\tau})\subset  \{x\in \t^*|\langle \xi, x\rangle \leq
\langle \xi, \mu \rangle\}.$$
In particular, if $a_w^{\tau}\neq 0$, then
\begin{equation}\label{eq:productlessthan}
\langle \xi, w\lambda\rangle \leq \langle \xi,\mu\rangle.
\end{equation}
Since $\xi$ is in $\tau\overline{C}$, it can be written as
\begin{equation}\label{eq:xi-in-tau-C}
\xi=\tau\sum_{j=1}^lr_j\lambda_j,
\end{equation}
where all $r_j$ are non-negative. So (\ref{eq:productlessthan}) and (\ref{eq:xi-in-tau-C}) imply that there exists
 $j\in \{1,\ldots,l\}$ such that
$$\langle \tau \lambda_j,w\lambda\rangle \leq \langle \tau \lambda_j,\mu\rangle.$$
In other words, each nonzero term in the right hand side of (3) is
a multiple of a $x_w^{\tau}$ of the type claimed in Theorem~\ref{th:main}. \qed

\noindent {\bf Remark.} It follows that, in the particular situation of generic
coadjoint orbits,  in order to cover the whole  Tolman-Weitsman kernel of the Kirwan
map  it is sufficient
to consider affine hyperplanes through $\mu$ which are perpendicular to vectors
of the type $\tau\lambda_j$, with $\tau\in W$ and $j\in \{1,\ldots, l\}$. But these
are just the hyperplanes parallel to the walls of the moment polytope.
This result concerning a ``sufficient set of hyperplanes''  has been proved by the
first author in \cite{Go2},
for an {\it arbitrary} Hamiltonian torus action on a compact manifold.

\bibliographystyle{abbrv}

\end{document}